
\input amstex

\def \restrict {\restriction}

\def \inj {\hookrightarrow}

\def\pp{\Bbb P}
\def\zz{\Bbb Z}
\def\cc{\Bbb C}
\def \qqed {\qed}

\catcode`\@=12
\def\cal{\Cal}
\def\SA{{\cal A}}

\def\SE{{\cal E}}
\def\SF{{\cal F}}
\def\SG{{\cal G}}

\def\SI{{\cal I}}

\def\SL{{\cal L}}
\def\SM{{\cal M}}
\def\SN{{\cal N}}
\def\SO{{\cal O}}

\def\SQ{{\cal Q}}

\input amstex

\documentstyle{amsppt}
\magnification=1200
\NoBlackBoxes

\rightheadtext\nofrills {Irreducibility of moduli of vector bundles}
\leftheadtext\nofrills {P.G.J. Nijsse}

\topmatter

\title
The irreducibility of the moduli space of stable vector bundles of rank
2 on a quintic in $\pp^3$
\endtitle

\author
Pieter G.J. Nijsse
\endauthor

\address Mathematical Institute Leiden University,
P.O.Box 9512,
NL-2300 RA Leiden,
The Netherlands \endaddress
\email nijsse\@wi.leidenuniv.nl \endemail

 \endkeywords

 \endsubjclass

\abstract {In this paper I consider a quintic surface in $\pp^3$, general
in the sense of Noether-Lefschetz theory. The vector bundles of rank 2
on this surface which are $\mu$-stable with respect to the hyperplane
section and have $c_1 = K$, the canonical class of the surface and
fixed $c_2$, are parametrized by a moduli space. This space is known to be
irreducible for large $c_2$ (work of K.G. O'Grady). I give an explicit
bound, namely $c_2 \geq 16$.} \endabstract

\endtopmatter

\def\mm{\SM (K, c_2)}
\def\nn{\SN (V)}

\document
\eject
\heading Introduction \endheading

Recently there has been some progress in proving irreducibility of the
moduli space of $\mu$-stable rank-2 vector bundles of large $c_2$, on a surface
of
general type. The first result in that direction was an ineffective
statement due to D. Gieseker and J. Li \cite {GL}, the second one was a
more effective result, giving a lower bound for $c_2$ to ascertain
irreducibility. This was proved by K.G. O'Grady in \cite {O2}. However, the
bound
that he obtained is only made effective for complete
intersections of high degree, and it is a hugh one. So the question
naturally arises whether one can do better in some simple cases.

In this paper I will consider quintic surfaces in $\pp^3$ general with respect
to the Noether-Lefschetz property, i.e. Pic $(S) = \zz$. For such a
surface it appears that the moduli space $\mm$ of $\mu$-stable (with respect to
the hyperplane section $H$) rank-2 vector bundles with Chern classes $c_1 = K$
and $c_2$, is irreducible for
$c_2 \geq 16$.
This is shown by establishing the appearance of non locally-free sheaves in the
closure of every component of $\mm$.

One may ask what happens in the case of a quintic surface $S$ with Pic $(S)
\not= \zz$.
O'Grady's proof (\cite {O2, 4.15-4.17}) of the irreducibility of $\mm$ does not
work here,
since it uses the existence of boundary points in the non-general case,
where we only have such a result for a general surface.
Although it seems not that difficult to redo the existence theorem,
allowing divisors that are not complete intersections
(the only place where it is really used occurs in step 3 to prove the
stability of the elementary modifications $\SF_y$), bounds will become
worse and I do not think it is worth the effort.
The main aim of this paper is to show that bounds on $c_2$ can be rather low in
special cases.

One also can ask what the situation is in case of surfaces of higher degree in
$\pp^3$, if one
looks at $\SM(H,c_2)$. This also can be attacked in the same way as
described in this paper, but again, bounds become worse, as could be
expected, whereas the proofs become more messy.

\vfill \eject

\heading \S 1 The boundary of the moduli space \endheading

First I need of course some notation. Let $S \subset \pp^3$ be a surface
of degree 5, general in the sense of Noether-Lefschetz, i.e. Pic $(S) =
\zz$. Let $K$ be the canonical class. Here $K$ is actually the
hyperplane class of $S$ and in view of this, it will sometimes be
denoted by $H$. Denote by $\mm$ the moduli space of rank-2 vector
bundles $\SE$ on $S$, with $c_1(\SE) = K$ and $c_2(\SE) = c_2$,
$\mu$-stable with respect to $H$. The existence and the construction of this
space can be found in \cite {M}.

Furthermore let $\overline \SM (K,c_2)$ be the closure of $\mm$ in the
Gieseker-Maruyama moduli space (see \cite {G})
of G-semistable torsion free sheaves. Define for a closed subset $X$ of the
moduli space $\overline \SM (K,c_2)$ the {\it boundary} of $X$ as

$$\partial X= \{ [\SF]
\in X \mid \SF \text { is not locally free} \}.$$
Notice that since $K$ is not divisible by 2, every $\mu$-semistable
sheaf is automatically $\mu$-stable.

I will take a 'component' to be an irreducible component. When I am
speaking about connected components I will always write 'connected
components'.

\proclaim {Proposition (1.1)}
Let $V= V(c_2) = \{ [\SE] \in \mm \mid H^0(\SE) \not= 0\}$.
 For $c_2 \geq 10$, $\dim V
\leq 3c_2 - 11 $.
\endproclaim

\demo {Proof}
If $\SE$ has a section, this gives rise to an exact sequence

$$ 0 \to \SO \to\SE \to \SI_Z(1) \to 0, $$

\noindent (see e.g. \cite {L, \S 3} for the well-known facts) where
$Z \in \text {Hilb}^{c_2}(S)$. Notice that I use that $\SE(-1)$ has no
sections for stability reasons, therefore the section of $\SE$ can not
vanish along a divisor.
Now let $\nn$ be the space

$$ \nn = \{ (\SE,s) \mid \SE\in V, s \in \pp H^0(\SE) \}$$
and consider the following diagram

$$
\CD
\nn @> p_2 >> \text {Hilb}^{c_2}(S) \\
@V p_1 VV @. \\
V @. \\
\endCD
$$
Certainly, $p_1$ is surjective. Hence dim $V \leq \dim \nn$. Furthermore

$$ \dim p_2^{-1}(Z) = \dim \pp \text {Ext}^1(\SI_Z(1), \SO) =
h^1(\SI_Z(2)) -1.$$
If $c_2 \geq 10$ a general $Z \in \text {Hilb}^{c_2}(S)$ has
$h^0(\SI_Z(2))=0$, so $h^1(\SI_Z(2)) = c_2 - 10$.
$\text {Hilb}^{c_2}(S)$ has dimension $2c_2$,
so the expected dimension of $\nn$ is $\leq 3c_2 -11$ (also valid for
$c_2=10$).

However we should be careful with the subsets

$$ \Delta_i = \{ Z \in \text {Hilb}^{c_2}(S) \mid h^0(\SI_Z(2)) \geq i \}.$$
In fact, their codimensions in $\text {Hilb}^{c_2}(S)$ should have to be $\geq
i$.
If we take inside
$$\pp H^0(\SO_S(2)) \times \text {Hilb}^{c_2}(S)$$
the incidence variety
$I = \{ (C,Z) \mid Z \subset C\}$ and denote by $\pi_1, \pi_2$ the projections,
then the dimension of $\pi_1^{-1} (C)$ is at most $c_2$, so
$9 + c_2 \geq \dim I \geq \dim \pi^{-1}_2(\Delta_i) \geq \dim \Delta_i + i-1$.
This implies that dim $\Delta_i$ is bounded by $c_2 + 10 - i \leq 2c_2 -i$
and hence the codimension of $\Delta_i$ in $\text {Hilb}^{c_2}(S)$ is $\geq i$.
\qqed
\enddemo

\proclaim {Theorem (1.2)} Let $t$ be a nonnegative integer.
Let $X \subset \overline \SM (K,c_2)$ be closed of dimension
$\geq 4c_2 - 20+t$. If $c_2 \geq \max
(10, 16 -t) $ then $\partial X
\not= \varnothing$.
\endproclaim

\demo {Proof}
We will follow O'Grady's approach \cite {O2} very closely. So let's
suppose that $\SE$ is locally free for all $[\SE] \in X$.

\proclaim {Step 1} For every smooth curve $C \in |H|$ there is a vector bundle
 $[\SE]
\in X$ such that $\SE \restrict_C$ is not $\mu$-stable, given that $c_2
\geq 9 - {1\over 4}t$.
\endproclaim

Let $\SM(C,H)$ be the moduli space of rank-2 semistable bundles on $C$
with determinant $H$. If $\SE \restrict_C$ were $\mu$-stable for all
$\SE$ then we would have a restriction map

$$ \rho : X \to \SM(C,H). $$

\noindent It is well known that the dimension of $\SM(C,H)$ is
$3(g_C-1)=15$. But for $c_2 \geq 9-{1\over 4}t$ the dimension of $X$ is at
least 16.
Hence, setting $\Theta$ the theta divisor
on $\SM(C,H)$ (see \cite {O2, prop.1.18}, \cite {DN}),

$$ (\rho^* \Theta)^{\text {dim} X} = 0,$$

\noindent contradicting the aforesaid proprosition 1.18 of \cite {O2}.
\bigskip
Set $X_C =\{ [\SE] \in X \mid \SE \restrict _C \text { is not stable}
\}$. This set is not empty as is shown in step 1. Recall that $V= \{ [\SE] \in
\mm \mid H^0(\SE) \not= 0\}$
as in proposition (1.1)

\proclaim {Step 2}
$X_C \backslash V \not= \varnothing$ provided that $c_2 \geq \max (10, 16 -t)$.
\endproclaim

Let $\SE$ be a vector bundle in $X_C$.
If $B$ is a neighbourhood of $[\SE]$ in $X$, the subset
$B^{\text {ns}} \subset B$ of vector bundles which are nonstable when
restricted to $C$ has codimension at most $g_C = 6$ by \cite {O2, prop.5.47}.
This means that $\dim B^{\text {ns}} \geq 4c_2 - 26+t$. But dim $V \leq 3c_2 -
11$.
So the claim
should hold for dimension reasons.

\proclaim {Step 3} The construction of a family of stable sheaves.
\endproclaim

Let $\SF$ be a vector bundle whose existence is guaranteed by step 2 and let

$$ 0 \to \SL_0 \to \SF \restrict_C \to \SQ_0 \to 0 $$

\noindent be a fixed destabilizing sequence. In particular  $c_1(\SL_0) >
c_1(\SQ_0)$
since $c_1(\SF \restrict _C)$ is not divisible by 2. To this sequence
elementary
 modifications of $\SF$ are
associated, namely, first consider the elementary transformation (\cite {L,
ex.3.17})

$$ 0 \to \SE \to \SF \to \iota_* \SQ_0 \to 0,$$

\noindent where $\iota : C \inj S$ is the natural injection. Notice
that we have $c_1(\SE) = c_1(\SF) - [C] = 0$ and
$H^0(\SE) \inj H^0(\SF)=0$, which implies
that $\SE$ is $\mu$-stable. Restricting this sequence to $C$ one gets

$$ 0 \to \SQ_0(-C) \to \SE \restriction _C \to \SL_0 \to 0.$$

Now let

$$ Y_{\SF} = Quot (\SE \restrict _C, \SL_0) $$
be the Grothendieck Quot-scheme parametrizing quotients of $\SE \restrict _C$,
that have the same Hilbert polynomial as $\SL_0$ (see
\cite {O2, p.8}). This set parametrizes a family
$\{\SF_y\}_{y \in Y_{\SF}}$ of
modifications of $\SF$, where $\SF_y = \SG_y (C)$ for a subsheaf $\SG_y$
of $\SE$, which is defined as the kernel in the exact sequence

$$ 0 \to \SG_y \to \SE \to \iota_*\SL_y \to 0. $$
Then it is easily shown that $\SF_y$ is $\mu$-stable
for all $y \in Y_{\SF}$, namely notice that
this is equivalent to the stability of $\SG_y$.  Then take a subsheaf $\SA=
\SO(a)
\inj \SG_y \inj \SE$. Since $\SG_y = \SF_y(-C)$, $\mu (\SG_y) = \mu (\SF_y) - 5
= - {5 \over 2}$. So we have to prove that $\mu
(\SA) < - {5 \over 2}$ i.e. $a < 0$. But this follows directly from the fact
that $\SE$ has no sections. Notice that this is the place where I make a
crucial use of the generality of $S$.

Thus we obtain a map

$$ \varphi : Y_{\SF} \to \overline \SM(K,c_2)$$

\noindent and therefore a subset $\varphi^{-1} (X) \subset Y_{\SF}$. The
dimension of this subset is bounded by the following elementary lemma

\proclaim {Lemma (1.3)}
$ \dim \varphi^{-1} (X) \geq \dim X + \dim Y_{\SF} - \dim T_{[\SF]} \mm.$
\endproclaim

\demo {Proof}
This is easily obtained by the following construction.
Notice that, since $\varphi$ is an injection, $\dim \varphi^{-1} (X) =
\dim X \cap Y_{\SF}$. Since $[\SF] \in  X \cap Y_{\SF}$, this space
is not empty. Embed a sufficiently small neighbourhood of $[\SF] \in \mm$ into
a
suitable $\cc^e$, namely with $e= \text {emb.dim}_{[\SF]} \mm$. By \cite {GR,
p.115} $e = \dim T_{[\SF]} \mm$. Then the lemma follows
from the result on intersection dimensions in affine space.
\qqed
\enddemo

The well known facts on $\mm$ (see e.g. \cite {OV} for a r\'esum\'e) as well
as the assumption that dim $X \geq \text {exp.dim}\ \mm$, imply that
$\dim T_{[\SF]} \mm - \dim X$ is bounded from above by
$h^2(\SE nd^0 \SF)= h^0(\SE nd^0 \SF(1))$.
Here $\SE nd^0 \SF$ denotes the traceless endomorphisms of $\SF$.
If we can arrange that the dimension of $\varphi^{-1} (X) > 1$ then \cite
{O2, lemma 1.15} assures the existence of a boundary point in $\partial
X$. Therefore lemma (1.3) shows that we are done if we prove the following
lemma

\proclaim {Lemma (1.4)}
$\dim Y_{\SF} > 1 + h^0(\SE nd^0 \SF(1))$
\endproclaim

\demo {Proof}
By definition of $Y_{\SF}$:

$$\split \dim Y_{\SF}
&= h^0(\SQ^* (C) \otimes \SL) \\
&= \chi (\SQ^* (C) \otimes \SL) + h^1(\SQ^* (C) \otimes \SL) \\
&= 1 + h^0(\SL^*\otimes \SQ(1)) \endsplit$$
(we skip the subscripts in $\SQ_0$ and $\SL_0$ here).
 So let's assume that

$$ h^0(\SE nd^0\SF(1)) > \dim \text {Hom} (\SL,\SQ(1)).$$
In view of the diagram

$$
\CD
0 @>>> \SL      @>\alpha >> {\SF \restrict _C}    @>>>  \SQ      @>>> 0  \\
@.     @.              @VsVV                     @.            @. \\
0 @>>> {\SL(1)} @>>>   {\SF(1) \restrict _C} @>\beta>> {\SQ(1)} @>>> 0  \\
\endCD
\tag 1
$$
this means that there is a section $s \in H^0(\SE nd^0 \SF(1))$ such
that $\beta \circ s \circ \alpha \equiv 0$, or equivalently that $s$ restricts
to
a map $a$ from $\SL$ to $\SL(1)$. Then $a$ can be viewed as an element of
$H^0(\SO_C(1))$ and since $C$ is a complete intersection, this comes
from an element of $H^0(\SO_S(1))$, which we will denote by $a$ too.

Consider $\tilde s = s -a: \SF \to \SF(1)$. Then $\SL \inj \ker \tilde s
\restrict _C$ and we have the following commutative diagram

$$
\CD
  @.                 @.       @.   0        @.      \\
@.     @.                 @.       @VVV          @. \\
  @.   0             @.       @.   \SO_D    @.      \\
@.     @VVV               @.       @VVV          @. \\
0 @>>> \SL           @>>> \SF \restrict _C @>>> \SQ      @>>> 0  \\
@.     @VVV               @|       @VVV          @. \\
0 @>>> \ker \tilde s \restrict _C @>>> \SF \restrict _C @>>> \text
{coker } \tilde s \restrict _C @>>> 0  \\
@.     @VVV               @.       @VVV            @. \\
  @.   \SO_D         @.       @.   0        @.      \\
@.     @VVV               @.       @.            @. \\
  @.   0             @.       @.            @.      \\
\endCD
$$
where $D$ is a divisor on $C$.
But since $\SQ$ is a destabilizing quotient of $\SF$ I can assume it to
be torsion free. Hence $D=0$, $\SL = \ker \tilde s \restrict _C$ and
the map $\SF \restrict _C \to \SQ$ is nothing else as $\tilde s \restrict
_C$. If we now take the elementary transformation associated to the composition
$ \SF \to \SF \restrict_C \to \SQ$ and denote it by $\SE$ then $\SE$ fits into
the exact sequence

$$ 0 \to \SE @>j>> \SF \to i_*\SQ \to 0. $$

\noindent So we obtain a map $\SE @>\tilde s \circ j >> \SE (1)$ which is zero
on $C$.
Hence either there is a non-trivial map $\SE(C) \to \SE(1)$ or $\tilde s \circ
j$ is the zero
map on $\SE$. In the first case one concludes that $ \tilde s \circ j \in
H^0(\SE nd\ \SE)$, which
contains only multiples of the identity by the stability of $\SE$.
In the second case $\tilde s$
is the zero map on $\SF$, i.e. $s = a$. In both cases there occurs a
contradiction with tr $s=0$.

Now it remains to exclude the case $h^0(\SE nd^0\SF(1)) = \dim \text {Hom}
(\SL,
\SQ(1))$. Diagram (1) shows that there is a map

$$ \psi : H^0 (\SE nd \ \SF (1)) \to \text {Hom} (\SL, \SQ(1)).$$
Moreover, the argument used above shows that $\psi$ restricted to the
subspace $H^0 (\SE nd^0 \SF (1))$ is injective. But if, as we have assumed,
the dimensions of these spaces match, the restricted $\psi$ is surjective too.
But that implies that the map $b: \SF \to \SF(1)$, where $b \in
H^0(\SO(1))$ has the same image under $\psi$ as a certain traceless map
$s_b$, say. This means that

$$\beta \circ (s_b -b) \circ \alpha \equiv 0.$$
 Now one can copy the above argument mutatis mutandis to obtain the
same contradiction.
\qqed \qqed
\enddemo
\enddemo

We close this section with the following proposition, showing that if there is
a boundary, it is large.

\proclaim {Proposition (1.5)}
(\cite {O2, prop.3.3}) Let $X \subset \overline \SM (K,c_2)$ be a closed
subspace
that has a non-empty boundary. Then $\partial X$ is closed and codim $(\partial
X, X) \leq 1$.
\qqed
\endproclaim

\vfill \eject

\heading \S 2 The goodness of the moduli space \endheading

We should start with saying what a good subset of the moduli space is

\proclaim {Definition (2.1)}
A closed subset $X \subset \overline \SM (K,c_2)$ is called {\it good} if in
every irreducible
component of $X$ there
is a sheaf $\SE$ with $H^2(\SE nd^0 \SE) =0$. Here $\SE nd^0 \SE$ denotes the
traceless endomorphisms of
$\SE$. {\it Bad} is not good.
\endproclaim

To simplify notation in the sequel I will use the following definition.

\proclaim {Definition (2.2)}
Let $X$ be a closed subset of $\overline \SM(K,c_2)$. Define

$$h(X) = \min \{h^2(\SE nd^0 \SF) \mid [\SF] \in X \}$$
\endproclaim

Notice that, if a component of the moduli space is good, it is reduced and of
the expected dimension
(exp.dim $\mm= 4c_2 -20$).
It turns out to be useful to know how large the vector space $H^2(\SE nd^0
\SF)$ can be,
for then we also know the maximal dimension of a component of $\mm$. By
Serre duality the dimension of this space equals $h^0 (\SE nd ^0 \SF
(1))$.

At a first glance one finds a bound from the sequence

$$ 0 \to \SO(-1) \to \SE nd^0 \SF \to \SG \to 0$$
for an appropriate $\SG$ (that one can prove to be a vector bundle).
Namely it implies that $h^2(\SE nd^0 \SF) \leq h^2(\SO(-1)) =10$ for all $\SF$.

However we can do much better when we have a vector bundle $\SF$, whose
restriction to a certain hyperplane section $C$ is not $\mu$-stable. The
existence of such an $\SF$ for large enough $c_2$ is guaranteed by step
1 of theorem (1.2).

\proclaim {Proposition (2.3)}
For $\SF$ a vector bundle, whose restriction to a certain hyperplane
section is not $\mu$-stable, $h^0(\SE nd^0\SF(1)) \leq 1$.
\endproclaim

\demo {Proof}
The proof of lemma (1.4) shows that $h^0(\SE nd^0 \SF(1)) < h^0(\SL^* \otimes
\SQ(1))$, where $\SL$ and $\SQ$ are the line bundles from the
destabilizing sequence on $C$. But the degree on $C$ of $\SL^* \otimes
\SQ(1)$ is less or equal than 4. Then Clifford's theorem \cite {ACGH,
p.107} implies that its space of sections has dimension $\leq 3$, with
equality if and only if $C$ is hyperelliptic. But $C$ is a complete
intersection, so it has ample canonical divisor and hence it is
not hyperelliptic. So $h^0(\SL^* \otimes
\SQ(1)) \leq 2$ and $h^0(\SE nd^0 \SF(1)) \leq 1$. \qqed
\enddemo

\proclaim {Corollary (2.4)}
If $X$ is a closed subset of $\mm$ of dimension at least 16 then $h(X)
\leq 1$.
\endproclaim

\demo {Proof}
Immediate from the proof of step 1 of theorem (1.2)
\qqed
\enddemo

Let me say some words about the {\it double-dual stratification} of $\partial
X$, for a closed subset $X$ of $\mm$.
For $[\SF] \in \partial X$ we have a canonical exact sequence

$$ 0 \to \SF \to \SF^{**} \to \SQ_{\SF} \to 0, $$

\noindent where $\SQ_{\SF}$ is a skyscraper sheaf of finite length $l
(\SQ_{\SF})=h^0(\SQ_{\SF})$.
The $\mu$-stability of $\SF$ implies
$\mu$-stability of $\SF^{**}$ (this is an easy lemma, see \cite
{Mo, lemma 2.2.1}). Notice that $c_2(\SF^{**}) = c_2(\SF) -
l(\SQ_{\SF})$.
Although in general it will not be possible to glue the bundles
$\SF^{**}$ for $\SF$ in $\partial X$ to a global family because their
second Chern classes may jump, there is a stratification of $\partial X$
by locally closed subsets, such that the duals of sheaves parametrized
by points of the same stratum, locally do glue to a neat family. This
follows from lemma 3.5 of \cite {O1} (note that one has locally
a universal sheaf on $\mm$). We are interested in
the open strata which have by proposition (1.5) codimension 1 in $X$.

\proclaim {Theorem (2.5)}
Let $t$ be a positive integer.
For $c_2 \geq \max (10, 16 -t)$, $\mm$ has no closed subsets $X$ of dimension
$\geq 4c_2 -20+t$ with $h(X) \geq 1$.
In particular, $\mm$ is good for $c_2 \geq 16$.
\endproclaim

\demo {Proof}
Let $X_0$ be a closed subset of $\mm$ of dimension $\geq 4c_2 -20+t$ and
 $h(X_0) \geq 1$. Since $c_2 \geq 9 - {1 \over 4} t$ by corollary (2.4)
$h(X_0) =1$ and hence $t \leq 1$ too. So we may assume that $c_2 \geq 15$.
Notice that $\partial X_0$ is not empty (theorem.(1.2)),
and proceed in the following way. Let
$Y_0$ be an irreducible component of an open stratum of the double-dual
stratification. Define $Y_0^{**} := \{\SF^{**} \mid [\SF] \in Y_0 \}$ and
set $X_1 = \overline {Y_0^{**}}$. Then $X_1$ lies in a certain moduli
 space $\SM(K,c'_2)$, with $c'_2 < c_2$. If still dim $X_1 \geq 16$ then
$h(X_1) = 1$ and by prop.3.8
together with cor.3.6 of \cite {O2} one gets

$$ \dim X_1 \geq 4c'_2 + (c_2 - c'_2) - 20 + t.$$

\noindent On the other hand dim $X_1 \leq 4c'_2 -20 +1$. So
$(c_2 -c'_2)+t \leq 1$, i.e. $c_2 -c'_2
=1$ and $t=0$. In particular $c'_2 \geq 15$ and dim $X_1 = 4c'_2 -20 +1$,
hence also $\partial X_1 \not= \varnothing$ and we are
able to construct $X_2$ in the given way. And again, if dim $X_2 \geq 16$,
then $(c_2 - c''_2) +t \leq 1$, which cannot occur.

So it remains to exclude the cases where dim $X_i < 16$ occurs for $i=1,2$.
Since by cor.3.6 of \cite {O2}

$$\split \dim X_i
& \geq 4c_2^{(i)} + (c_2 - c_2^{(i)}) -21 +t \\
& \geq 3c_2^{(i)} -5,\endsplit
$$
this corresponds to $c_2^{(i)} < 7$. But in that case dim $X_i \leq \dim V \leq
3c_2^{(i)} -11$ by proposition (1.1).
\qqed
\enddemo

\vfill \eject

\heading \S 3 The irreducibility of the moduli space \endheading

Let $X$ and $Y$ be two irreducible components of $\mm$,
meeting in a subset $Z$. Assume that both $X$ and $Y$ are good and have
boundaries.
If their boundaries do not meet each other they are coming from
components $X_1$ and $Y_1$ of $\SM(K,c_2-1)$, whose closures are disjoint.
Thus $\overline \SM (K,c_2-1)$ will be disconnected. If on the contrary,
both boundaries do meet, $Z$ has a boundary, which gives rise to a bad
component of $\SM (K,c_2-1)$. This shows that to prove
irreducibility it will suffice to prove that $\overline \SM (K,c_2-1)$
is connected. Therefore it is useful to have the following lemma

\proclaim {Lemma (3.1)}
$\overline V(c_2) := \{[\SF] \in \overline \SM(K,c_2) \mid h^0(\SF) \not= 0 \}$
is connected
for $c_2 \geq 10$.
\endproclaim

\demo {Proof}
If $\SF$ has a section, the cokernel of that section is torsion free
since $\mu (\SO) = 0$ is maximal in the set $\{\mu (\SA) \mid \SA \inj
\SF \}$. So exactly as in the case of a vector bundle one obtains the
exact sequence

$$ 0 \to \SO \to \SF \to \SI_Z(1) \to 0 $$
and we can define the same sets and maps as we did in proposition (1.1).
Notice that this uses the generality in the sense of Noether-Lefschetz
of $S$, however we will see at the end of this section that this
dependence is not essential.

For $c_2 \geq 11$, a general point in $\text {Hilb}^{c_2}(S)$ has the
Cayley-Bacharach property with respect to quadrics, so $p_2(\SN(\overline V))
\subset
\text {Hilb}^{c_2}(S)$ is dense and hence connected. Furthermore, the
fibres of $p_2$ are the projective spaces $\pp \text {Ext}^1(\SI_Z(1),
\SO)$.

For $c_2 =10$, $\SN (\overline V)$ is mapped by $p_2$ onto the set $\{ Z \mid
h^0(\SI_Z(2)) \not= 0 \} \subset \text {Hilb}^{10} (S)$, with connected
fibres. But this set is certainly connected.
\qqed
\enddemo

\proclaim {Proposition (3.2)}
For $c_2 \geq 10$, $\overline \SM (K,c_2)$ is connected.
\endproclaim

\demo {Proof}
Notice that the proof of theorem (1.2) shows that an irreducible component $X$
of $\mm$, for $c_2 \geq 10$, has a boundary unless $V(c_2) \cap X \not=
\varnothing$. So the moduli space $\overline \SM(K,c_2)$ has a connected
component $C(V)$ containing $\overline V$
and all other connected components have boundaries. The latter are
coming from a moduli space $\overline \SM(K,c'_2)$, with $c'_2 = c_2-1$.
But for $c_2 \geq 11$ $C(V)$ has a non-empty boundary too, since for $Z$ being
$c_2 -1$
points on a quadric plus a point not on that quadric, certainly
$h^1(\SI_Z(2)) \not= 0$, but $Z$ does not have the Cayley-Bacharach
property with respect to quadrics, so the torsion free sheaf defined by
$Z$ is not locally free.
So the number of connected components of $\overline \SM(K,c_2)$ is a decreasing
function of $c_2$ for $c_2 \geq 10$.

So it remains to show that $\overline \SM(K,10)$ is connected. As is
said already, we certainly have the connected component $C(V)$. Assume that
there is another connected component $X$.
Then $X$ has a boundary and we can choose an element $[\SF] \in \partial X$.
$\SF$
gives rise to $\SF^{**}$, which fits into

$$ 0 \to \SO \to \SF^{**} \to \SI_Z(1) \to 0,$$
since by Riemann-Roch $h^0(\SF^{**}) \geq 10 - c_2(\SF^{**}) >0$. Now
choose a set $W$ of $l(\SQ_{\SF})$ different points, disjoint from $Z$. Let
$\SG$ be the torsion free sheaf corresponding to the image of
the extension class of $\SF^{**}$ under the map

$$\text {Ext}^1(\SI_Z(1), \SO) \to \text {Ext}^1 (\SI_{Z+W}(1), \SO).$$
This means that $\SG \inj \SF^{**}$ is defined as the kernel of the
composite map

$$\SF^{**} \to \SI_Z(1) \to \SO_W.$$

Thus we obtain the following commutative diagram

$$ \CD
@.@. 0 @. 0 @.\\
@.@. @VVV @VVV @.\\
0 @>>> \SO @>>> \SG @>>> \SI_{Z+W}(1) @>>> 0\\
@. @| @VVV @VVV @.\\
0 @>>> \SO @>>> \SF^{**} @>>> \SI_Z(1) @>>> 0 \\
@. @. @VVV @VVV @. \\
@. @. \SO_W @= \SO_W @.\\
@.@. @VVV @VVV @.\\
@. @. 0 @. 0 @. \\
\endCD $$
which shows that $\SG$ sits in the same connected component of the
boundary $\partial X$ as $\SF$, since $W$ can be obtained from a
continuous disturbation of the points in the support of $\SQ_{\SF}$.
Moreover it shows that $H^0(\SG) \not= 0$. This contradicts our assumption
that $X$ is a connected component different from $C(V)$.
\qqed
\enddemo

But by what is said above we can make an improvement of theorem (2.6).

\proclaim {Theorem (3.3)}
$\mm$ is good for $c_2 \geq 13$
\endproclaim

\demo {Proof}
Let $X$ be a bad component and assume that there is another component
$Y$. Since the moduli space is connected, they intersect. But $h(X \cap
Y)$ will be greater than $h(X)=1$, so by corollary (2.4) dim $X \cap Y
\leq 15$, or equivalently, its codimension in $Y$ is $\geq 4c_2 - 35$.
But this codimension may not exceed 10, because $h^0(\SE nd^0 \SF(1))$
is bounded by 10. So $c_2 \leq 11$, a contradiction.

This shows that $X$ is the only component of $\mm$. But then it
certainly has a boundary, coming from a bad component $X_1$ of $\SM(K,c_2-1)$
and by the same argument $X_1$ is coming from an even worse component
$X_2$, which does not exist (see the argument of theorem (2.5)).
\qqed
\enddemo

The final and main theorem is

\proclaim {Theorem (3.4)}
$\mm$ is irreducible for $c_2 \geq 16$.
\endproclaim

\demo {Proof}
This is easy now: it is just the remark at the beginning of this
section, using theorem (3.3).
\qqed
\enddemo

\proclaim {Remark (3.5)} The situation that can a priori occur for
$c_2 =15$ is that there are many good components without boundary.
\endproclaim

Now the question arises what we can do in case of a surface $S$ with Pic
$S \not= \zz$.
We have used that Pic $S=\zz$ in lemma (3.1), where I indicated that
proposition (3.2) can be proven without this assumption. Indeed, lemma
(3.1) shows that in any case the closure of $V(c_2)$ is connected. So
when Pic $S \not= \zz$, it remains possible that there are boundaries,
disjoint from $C(V)$, whose corresponding sheaves have sections. But the
diagram at the end of the proof of proposition (3.2) shows that after
some small disturbation the cokernel of the section is torsion free. So
it has to meet $C(V)$ and proposition (3.2) remains valid even in the
non-general case.

As already mentioned in the introduction the only place
where the assumption $\text {Pic} S = \zz$ further is used, is in step 3 of
theorem
(1.2), where it is indicated. So it is used to find boundary points in a
component. However the proof of theorem (3.3) shows that you do not need
to prove the existence of boundary points in case you want to prove that
the moduli space is good, since in bad components you obtain your
boundary points for free. So theorem (3.3) holds for an arbitrary quintic
surface. However this is not enough to prove the irreducibility. I
really need boundary points there (see \cite {O2, 4.15-4.17}).

\enddocument